\documentclass[10pt]{article}
\usepackage{amsmath,amsthm,amscd,amssymb}
\usepackage{latexsym}

\newtheorem{theorem}{Theorem}[section]

\newtheorem{corollary}[theorem]{Corollary}

\theoremstyle{definition}

\theoremstyle{remark}

\numberwithin{equation}{section}
\allowdisplaybreaks

 \makeatletter
\begin{document}

\title{ Applications of Some Elliptic Equations in Riemannian Manifolds
       \thanks{This work  was partly supported by Fujian Scholarship Fundation and NFS of Fujian Province (Grant No. 2012J01015).} }

\author{{\small Qin Huang   \thanks { \textit{E-mail address}: qinhuang78@163.com   }
         \ \ \ Qihua Ruan\thanks{\textit{Corresponding author}: ruanqihua@163.com}
} \\
%EndAName
{\small  Department of Mathematics,
Putian University, Putian, 351100, Fujian, PR China}\date
{\noindent\small }\\
}
\date{}
\maketitle
%\address{3}%
%\thanks{4}%
\date{}%
% ----------------------------------------------------------------
\begin{abstract}
 Let $(M^{n+1}, g)$ be a compact Riemannian manifold with smooth boundary B and nonnegative Bakry-Emery Ricci curvature. In this paper, we use the solvability of some elliptic equations to prove some estimates of the weighted mean curvature and some related rigidity theorems. As their applications, we obtain some  lower bound estimate of the first nonzero eigenvalue of the drifting Laplacian acting on functions on B and some corresponding rigidity theorems.



\end{abstract}
%\maketitle
% ----------------------------------------------------------------
\section{Introduction and main results}
Let $(M^{n+1}, g)$ be a compact Riemannian manifold with smooth boundary $B^{n}=\partial M$.
Let $dV$ and $dA$ be the canonical measures on M and B respectively ,  V and A be the volume of
M and the area of B. Given $f\in C^{\infty}(M)$, Reilly's formula \cite{R2} states that
\begin{eqnarray}
\int_{M}((\bar{\triangle}f)^{2}-|\bar{\nabla}^{2}f|^{2}-Ric(\bar{\nabla }f, \bar{\nabla }f))dV
=\int_{B}(2(\triangle f)f_{\nu}+nH(f_{\nu})^{2}+\Pi(\nabla f, \nabla f))dA,\label{1.1}
\end{eqnarray}
here $\bar{\nabla}f$, $\bar{\triangle}f$, $\bar{\nabla}^{2}f$ being the gradient, the Laplacian and
the Hessian of $f$ on M, $Ric$  the Ricci curvature of M, $\nabla f$ and $\triangle f$ the gradient
and the Laplacian of $f$ in B, and $\Pi(X, Y)=g(\bar{\nabla}_{X}\nu, Y)$ for $\forall X, Y\in TB$ and $H=\frac{1}{n}tr \Pi$
the second fundamental form and the mean curvature
of B with respect to the outer unit normal $\nu$ on B.

%Let $f:M\rightarrow R$ be a solution of the following elliptic equation with Neumann boundary

%\begin{equation*}\bar{\triangle} f=1   \ \ \text{on} \  M, \  f_{\nu}(x)=\frac{V}{A}  \ \  \text{on}\  B. \end{equation*}
 Reilly \cite{Reilly1} used the formula \eqref{1.1} to prove that if M has non-negative Ricci curvature with  convex boundary
  and $H\geq\frac{A}{(n+1)V}$, then M is isometric to an Euclidean ball. Later Ros \cite{Ros} removed the condition of the convex boundary
  and obtains the same conclusion.

Recently Ma and Du \cite{MD} studied the drifting laplacian operator $\bar{\Delta}_{h}=\bar{\triangle}-\bar{\nabla}h\cdot\bar{\nabla}$ for a
smooth function $h$ on M. This operator is self-adjoint operator with respect to the weighted measure $dV_{h}=e^{-h}dV$. They extended the above
Reilly's formula and showed that
\begin{eqnarray}
\int_{M}((\bar{\triangle}_{h}f)^{2}-|\bar{\nabla}^{2}f|^{2}-Ric_{h}(\bar{\nabla }f, \bar{\nabla }f))dV_{h}
=\int_{B}(2(\triangle_{h} f)f_{\nu}+nH_{h}(f_{\nu})^{2}+\Pi(\nabla f, \nabla f))dA_{h}\label{1.2}
\end{eqnarray}
here $Ric_{h}=Ric+\bar{\nabla}^{2}h$, $\triangle_{h}f=\triangle f -\nabla h\cdot\nabla f$,  $H_{h}=H-\frac{1}{n}h_{\nu}$
and $dA_{h}=e^{-h}dA.$ In \cite{Ruan2}, the second author proved the sharp gradient estimate for positive solution of the Laplacian with a general dift $B$, i.e.
$\bar{\triangle}f-Bf=0,$ where B is a vector field.

Using the inequalities of $|\bar{\nabla}^{2}f|^{2}\geq\frac{1}{n+1}(\bar{\triangle}f)^{2}$ and $\frac{1}{n+1}a^{2}+\frac{1}{m-n-1}b^{2}
\geq\frac{1}{m}(a-b)^{2}$, we know that the equation \eqref{1.2} become the following inequality:
\begin{eqnarray}\int_{M}\frac{m-1}{m}((\bar{\triangle}_{h}f)^{2}-Ric_{m}(\bar{\nabla }f, \bar{\nabla }f))dV_{h}
\geq\int_{B}(2(\triangle_{h} f)f_{\nu}+nH_{h}(f_{\nu})^{2}+\Pi(\nabla f, \nabla f))dA_{h}\label{1.3}
\end{eqnarray}
here $Ric_{m}=Ric_{h}-\frac{1}{m-n-1}\bar{\nabla }h\bigotimes\bar{\nabla }h, m\geq n+1,$ and $m=n+1$ if and only if $h$ is a constant.  This curvature tensor is called Bakry-Emery Ricci curvature
(see \cite{Li}).
 With help of this inequality, Ma and Du obtained the lower
bound for the first eigenvalue of the drifting Laplacian on the compact manifold with positive Bakry-Emery Ricci curvature( see \cite{MD}). It states that if $Ric_{m}\geq(m-1)K>0$, $H_{h}\geq0$ or $\Pi\geq 0$, then    the first Neumann eigenvalue $\lambda^{N}_{1}(\bar{\triangle}_{h})\geq mK$, or the first Dirichlet eigenvalue $\lambda^{D}_{1}(\bar{\triangle}_{h})\geq mK.$ This conclusion is a generalization of Reilly's \cite{Reilly1} and Escobar's results \cite{E}. Recently Li and Wei \cite{LW} proved that this result is sharp, i.e. if $\lambda^{D}_{1}(\bar{\triangle}_{h})= mK$ or $\lambda^{N}_{1}(\bar{\triangle}_{h})= mK$, then the manifold is isometric to  a Euclidean hemisphere. They extended the rigidity theorem of Reilly\cite{Reilly1} and Escobar\cite{E}.  From  Theorem 5 in \cite{Q}, we know that if $Ric_{m}\geq(m-1)K>0$, then M is compact and the diameter $diam(M)\leq\frac{\pi}{\sqrt{K}}$. In \cite{Ruan1}, the second author proved that if $Ric_{m}\geq(m-1)K>0$ and $diam(M)=\frac{\pi}{\sqrt{K}}$, then M is isometric to a Euclidean sphere of radius $\frac{1}{\sqrt{K}}.$

Based on the Reilly formula \eqref{1.1}, Ros \cite{Ros} showed an estimate of the mean curvature. Similarly using the Reilly type inequality \eqref{1.3} we may  extend Ros's  result to the manifold with nonnegative Bakry-Emery Ricci curvature and obtain an estimate of the weighted mean curvature. However in this paper we do not use  the Reilly type formula but the divergence theorem   to prove this result.

%theorem 1.1
\begin{theorem}Let $(M^{n+1}, g)$ be a compact Riemannian manifold with smooth boundary B and $Ric_{m}\geq0$. If  the weighted mean curvature $H_{h}>0$, then
\begin{eqnarray}
\int_{B}\frac{dA_{h}}{H_{h}}\geq\frac{mn}{m-1}V_{h},\label{1.4}
 \end{eqnarray}here $V_{h}=\int_{M}dV_{h}$ denotes the weighted volume of M. The equality holds if and only if M is isometric to an Euclidean ball and $h$ is constant.
\end{theorem}
Let $A_{h}=\int_{B}dA_{h}$,  if $H_{h}\geq\frac{(m-1)A_{h}}{mn V_{h}}$, then the equality in \eqref{1.4} holds. Thus we  easily deduce the following rigidity theorem.
%corollary 1.2
\begin{corollary}Let $(M^{n+1}, g)$ be a compact Riemannian manifold with smooth boundary B and $Ric_{m}\geq0$. If  the weighted mean curvature $H_{h}\geq\frac{(m-1)A_{h}}{mn V_{h}}$, then M is isometric to an Euclidean ball and $h$ is constant.

\end{corollary}
{\bf Remark 1:} When $h$ is a constant, the Bakry-Emery Ricci curvature and the weighted mean curvature become the classical Ricci curvature and the mean
curvature respectively, and $m=n+1$. In this case, Corollary 1.2 is Ros's result\cite{Ros}.

Using the similar method we  prove the following estimate of the weighted mean curvature and the related rigidity theorem. If $h$ is a constant, then
it is a  Reilly's result in \cite{Reilly1}.
%theorem 1.3
\begin{theorem}Let $(M^{n+1}, g)$ be a compact Riemannian manifold with convex boundary B and $Ric_{m}\geq0$. Then \begin{eqnarray}\int_{B}H_{h}dA_{h}\leq\frac{(m-1)A^{2}_{h}}{mn V_{h}},\label{1.5}
\end{eqnarray}  The equality holds if and only if M is isometric to an Euclidean ball and $h$ is constant.
\end{theorem}

Now we discuss some applications of Corollary 1.2. Firstly we prove that if the manifold has nonnegative Bakry-Emery
 Ricci curvature, then the critical point of the weighted isoperimetric functional is an Euclidean ball. Recently there are many results about
 isoperimetric problems on the manifold with density, for example, see  \cite{CMCR}, \cite{Espinar}, \cite{Ho}, \cite{MM},\cite{Morgan}, \cite{RCBM}  and so on.
%theorem 1.4
\begin{theorem}Let $(M^{n+1}, g)$ be a compact Riemannian manifold with $Ric_{m}\geq0$, and let $\Omega$ be a compact domain in M with smooth boundary $\partial\Omega$. If $\Omega$ is a critical point of the weighted isoperimetric functional
 $$\Omega\rightarrow \frac{A_{h}(\partial \Omega)^{m}}{V_{h}(\Omega)^{m-1}},$$ then $\Omega$ is isometric to an Euclidean ball and $h$ is constant.
\end{theorem}

In \cite{xia}, Xia   used the Reilly's formula \eqref{1.1} and Ros's result\cite{Ros} to obtain a lower bound of the first nonzero eigenvalue $\lambda_{1}(\triangle)$ of the Laplacian acting on
functions on B and the corresponding  rigidity theorem.
%theorem 1.5
\begin{theorem}{\bf (Xia's theorem in \cite{xia})} Let $(M^{n+1}, g)$ be a compact Riemannian manifold with nonempty boundary B and nonnegative Ricci curvature. If  the second fundamental form of B satisfies $\Pi\geq cI$ (in the matrix sense), then  $\lambda_{1}(\triangle)\geq nc^{2}$. The equality holds if and only if M is isometric to an Euclidean ball.
\end{theorem}

In this paper, we use different method to generalize this result to the drifting Laplacian $\triangle_{h}=\triangle-\nabla h\cdot\nabla$. Our proof
is mainly based on the divergence theorem  and Corollary 1.2.
%theorem 1.6
\begin{theorem}Let $(M^{n+1}, g)$ be a compact Riemannian manifold with nonempty boundary B and $Ric_{m}\geq0$. If  the second fundamental form of B satisfies $\Pi\geq cI$ (in the matrix sense) and $h_{\nu} \leq -(m-n-1)c$, then $\lambda_{1}(\triangle_{h})\geq (m-1)c^{2}$. The equality holds if and only if M is isometric to an Euclidean ball of radius $\frac{1}{c}$ and $h$ is constant.
\end{theorem}

Using the similar method in the proof of Theorem 1.6, we also show  the  following result.
%theorem 1.7
\begin{theorem}Let $(M^{n+1}, g)$ be a compact Riemannian manifold with nonempty boundary B and $Ric_{m}\geq0$. If  the second fundamental form of B satisfies $\Pi\geq cI$ (in the matrix sense), and the weighted mean curvature $H_{h}\geq\frac{\lambda_{1}(\triangle_{h})}{nc}$, then $\lambda_{1}(\triangle_{h})\leq(m-1)c^{2}$.
\end{theorem}
When $h$ is a constant, we obtain the following corollary.
%corollary 1.8
\begin{corollary}Let $(M^{n+1}, g)$ be a compact Riemannian manifold with nonempty boundary B and nonnegative Ricci curvature. If  the second fundamental form of B satisfies $\Pi\geq cI$ (in the matrix sense) and the  mean curvature $H\geq\frac{\lambda_{1}(\triangle)}{nc}$, then $\lambda_{1}(\triangle)\leq nc^{2}$.
\end{corollary}
Combining Corollary 1.8 with Theorem 1.5, we know that under the condition of Corollary 1.8 the manifold M is isometric to an Euclidean ball .
 This  is the Theorem 3 in \cite{xia} proved by Xia.
%corollary 1.9
\begin{corollary}{\bf (Xia's theorem in \cite{xia})}Let $(M^{n+1}, g)$ be a compact Riemannian manifold with nonempty boundary B and nonnegative Ricci curvature. If  the second fundamental form of B satisfies $\Pi\geq cI$ (in the matrix sense) and the  mean curvature $H\geq\frac{\lambda_{1}(\triangle)}{nc}$, then M is isometric to an Euclidean ball of radius $\frac{1}{c}$ .
\end{corollary}

{\bf Acknowledgement.} The authors wish to express their thanks to
Professor X.N. Ma and the anonymous referee for pointing out some errors in the first manuscript and giving us some important suggestions. The authors also thank professor J.M.Espinar for sending us his new version  of the reference \cite{Espinar}.

\section{Proof of  Theorem 1.1, 1.3 and 1.4}
{\bf Proof of Theorem 1.1.} Let $f$ be the smooth solution of the Dirichlet problem

\begin{equation*}  \left \{
\begin{array}{lll}
\bar{\triangle}_{h} f&=1 &\text{ in }  M,   \\
f&=0  &\text{ on } B.
\end{array}
\right.
\end{equation*}
Suppose $F=\frac{1}{2}|\bar{\nabla }f|^{2}-\frac{1}{m}f$, then we have that
\begin{eqnarray}{\label{2.1}}
  \bar{\triangle}_{h}F&=& |\bar{\nabla}^{2}f|^{2}+\bar{\nabla}\bar{\triangle}_{h}f\bar{\nabla}f+Ric_{h}(\bar{\nabla }f,\bar{\nabla }f)-\frac{1}{m}\bar{\triangle}_{h}f \\
  &\geq&\frac{1}{n+1}(\bar{\triangle}f)^{2}+\frac{1}{m-n-1}(\bar{\nabla}h\bar{\nabla}f)^{2}+\bar{\nabla}\bar{\triangle}_{h}f\bar{\nabla}f+Ric_{m}(\bar{\nabla }f,\bar{\nabla }f)-\frac{1}{m}\bar{\triangle}_{h}f \nonumber \\
  &\geq& \frac{1}{m}(\bar{\triangle}_{h}f)^{2}+\bar{\nabla}\bar{\triangle}_{h}f\bar{\nabla}f+Ric_{m}(\bar{\nabla }f,\bar{\nabla }f)-\frac{1}{m}\bar{\triangle}_{h}f\nonumber \\
 &\geq& \frac{1}{m}(\bar{\triangle}_{h}f)^{2}+Ric_{m}(\bar{\nabla }f,\bar{\nabla }f)-\frac{1}{m}\bar{\triangle}_{h}f\geq0.\nonumber
\end{eqnarray}

Integrating both sides of the above inequality on M with respect to the weighted  measure $dV_{h}$,  we have that
$$\int_{M}\bar{\triangle}_{h}F(x)dV_{h}\geq0.$$

By the divergence theorem ,  we know that $$\int_{B}\frac{\partial F}{\partial\nu}(x)dA_{h}\geq0.$$ Since $f=0$ on B, then $\nu=\frac{\bar{\nabla}f}{|\bar{\nabla}f|}$ and $F(x)=\frac{1}{2}f_{\nu}^{2}-\frac{1}{m}f.$ From the fact that $\triangle_{h} f+f_{\nu\nu}+nH_{h}f_{\nu}=1$ and $f=0$ on B, we conclude that
\begin{equation*}
0\leq \int_{B}F_{\nu}dA_{h}=\int_{B}(f_{\nu}f_{\nu\nu}-\frac{1}{m}f_{\nu})dA_{h}=\int_{B}(\frac{m-1}{m}f_{\nu}-nH_{h}f_{\nu}^{2})dA_{h}.
\end{equation*}
Thus we have that
\begin{eqnarray}
\int_{B}H_{h}f_{\nu}^{2}dA_{h}\leq\frac{m-1}{mn}V_{h}. \label{2.2}
\end{eqnarray}
Here we use the following equation.
\begin{eqnarray}
\int_{B}f_{\nu}dA_{h}=\int_{M}\bar{\triangle}_{h}fdV_{h}=V_{h}. \label{2.3}
\end{eqnarray}
Finally, from \eqref{2.2}, \eqref{2.3} and Schwarz inequality it follows that
\begin{eqnarray*}
  V_{h}^{2}&=&\left(\int_{B}f_{\nu}dA_{h}\right)^{2}=(\int_{B}(H_{h}^{\frac{1}{2}}f_{\nu})H_{h}^{-\frac{1}{2}}dA_{h})^{2} \\
  &\leq&\int_{B}H_{h}f_{\nu}^{2}dA_{h}\int_{B}H^{-1}_{h}dA_{h}\leq\frac{m-1}{mn}V_{h}\int_{B}H^{-1}_{h}dA_{h}.
\end{eqnarray*}

Thus we have proved the inequality \eqref{1.4}.

If M is isometric to an Euclidean ball and $h$ is constant, then it is easy to conclude that the equality sign in (1.4) holds. Now we assume conversely that the equality sign in \eqref{1.4} holds. In this case all the equalities hold in \eqref{2.1}. Thus $\frac{\bar{\triangle}f}{n+1}=-\frac{\bar{\nabla}h\bar{\nabla}f}{m-n-1}=\frac{\bar{\triangle}_{h}f}{m}$. As the proof of Theorem 3 in \cite{LW}, we know that $m=n+1$, $h$ is constant. In fact, we have the following equation
 \begin{eqnarray}
\bar{\triangle}f=-\frac{n+1}{m-n-1}\bar{\nabla}h\bar{\nabla}f.\label{2.4}
 \end{eqnarray}
 If $m>n+1,$ then multiplying \eqref{2.4} with $f$ and integrating on M with respect to $e^{\frac{n+1}{m-n-1}h}dV$, we obtain that $\int_{M}|\bar{\nabla}f|^{2}e^{\frac{n+1}{m-n-1}h}dV=0.$ So $f$ is constant, which is a contradiction since $\bar{\triangle}_{h}f=1.$ Thus by Ros's result we see that M is isometric to an Euclidean ball.

{\bf Proof of Theorem 1.3.} Let $f$ be the smooth solution of the Neumann problem

\begin{equation*}  \left \{
\begin{array}{lll}
\bar{\triangle}_{h} f&=1 &\text{ in }  M,   \\
f_{\nu}&=\frac{V_{h}}{A_{h}}  &\text{ on } B.
\end{array}
\right.
\end{equation*}
Suppose $F=\frac{1}{2}|\bar{\nabla }f|^{2}-\frac{1}{m}f$,  from the proof of Theorem 1.1, we know that
\begin{eqnarray}\bar{\triangle}_{h}F\geq0.\label{2.5}
\end{eqnarray} Then by the divergence theorem,  we know that $\int_{B}\frac{\partial F}{\partial\nu}(x)dA_{h}\geq0.$ On the other hand,
\begin{eqnarray} \label{2.6}
  F_{\nu}&=&\sum\limits_{i=1}^{n+1} f_{i\nu}f_{i}-\frac{1}{m}f_{\nu}\nonumber\\
  &=& \sum\limits_{i=1}^{n} f_{\nu i}f_{i}-\Pi_{ij}f_{i}f_{j}+f_{\nu\nu}f_{\nu}-\frac{1}{m}f_{\nu} \nonumber \\
  &\leq& f_{\nu\nu}f_{\nu}-\frac{1}{m}f_{\nu} \\
  &=&(1-\triangle_{h}f-nH_{h}f_{\nu})f_{\nu}-\frac{1}{m}f_{\nu}\nonumber \\
  &=&\frac{m-1}{m}f_{\nu}-\triangle_{h}ff_{\nu}-nH_{h}f_{\nu}^{2}.\nonumber
\end{eqnarray}

Thus from the boundary condition and the divergence theorem,  we have that
\begin{eqnarray}\int_{B}H_{h}dA_{h}\leq \frac{m-1}{nm}\frac{A_{h}^{2}}{V_{h}}.
\end{eqnarray}

If M is isometric to an Euclidean ball and $h$ is constant, then it is easy to conclude that the equality sign in \eqref{1.5} holds. Now we assume conversely that the equality sign in \eqref{1.5} holds. As the proof of Theorem 1.1, we obtain that $\frac{\bar{\triangle}f}{n+1}=-\frac{\bar{\nabla}h\bar{\nabla}f}{m-n-1}=\frac{\bar{\triangle}_{h}f}{m}$. We claim that $m=n+1$, $h$ is constant. In fact, we have the following equation $$\bar{\triangle}f=-\frac{n+1}{m-n-1}\bar{\nabla}h\bar{\nabla}f.$$ If $m>n+1,$ then  integrating the above equation on M with respect to $e^{\frac{n+1}{m-n-1}h}dV$, we obtain that $\int_{B} f_{\nu}e^{\frac{n+1}{m-n-1}h}dA=0$. This is a contradiction since $f_{\nu}=\frac{V_{h}}{A_{h}}.$ On the other hand, if the equality sign in \eqref{1.5} holds, then  the equalities in \eqref{2.5} and \eqref{2.6} hold. Which implies that $\bar{\nabla}^{2}f=\frac{\bar{\triangle}f}{n+1}g$ and $\Pi(\bar{\nabla}f, \bar{\nabla}f)=0$, i.e.
\begin{eqnarray}\bar{\nabla}^{2}f&=&\frac{1}{n+1}g \ \ \ \text{in} \ M; \nonumber\\f&=&\text{constant}\ \ \  \text{on} \  B.\nonumber
\end{eqnarray}
Then by Lemma 3 in \cite{Reilly1} we see that M is isometric to an Euclidean ball.

The proof of Theorem 1.4 is almost following Ros's method \cite{Ros} and the first variation formulae of weighted volume and perimeter \cite{RCBM}. We include here a proof for the sake of completeness.

{\bf Proof of Theorem 1.4. } Given a smooth function $u$ on $\partial \Omega$, we consider the normal variation of $\partial \Omega$ defined by $\varphi_{t}: \partial \Omega\rightarrow M$, $\varphi_{t}(p)=Exp_{p}(tu(p)N(p)),$ where $Exp$ is the exponential map of M. $\varphi_{t}$ determine a variation of $\Omega$, $\Omega_{t}$ for $|t|<\epsilon$. Let $V_{h}(t)=V_{h}(\Omega_{t})$ and $A_{h}(t)=A_{h}(\partial\Omega_{t})$. The first variation formulae of the functionals above are given by
\begin{eqnarray*}
A_{h}'(0)&=&n\int_{\partial\Omega}uH_{h}dA_{h} \\
V_{h}'(0)&=&\int_{\partial\Omega}udA_{h}.
\end{eqnarray*}

By hypothesis we show
$$\frac{d}{dt}\mid_{t=0}\frac{A_{h}(t)^{m }}{V_{h}(t)^{m-1}}=0,$$
or equivalently
\begin{equation*}\int_{\partial\Omega}u((mnV_{h}H_{h}-(m-1)A_{h}))dA_{h}=0, \text{ for any} \ u.
\end{equation*}

 Then $H_{h}=\frac{(m-1)A_{h}}{mnV_{h}}$. Therefore from Corollary 1.2 $\Omega$ is isometric to an Euclidean ball and $h$ is constant.

\section{Proof of  Theorem 1.6, 1.7 and 1.11 }

{\bf Proof of Theorem 1.6. } Let $u$ be an eigenfunction corresponding to the first nonzero eigenvalue $\lambda_{1}$ of the drifting Laplacian of B:
\begin{equation*}\triangle_{h}u=-\lambda_{1}u.
\end{equation*}
Let $f$ be the smooth solution of the Dirichlet problem:
\begin{equation*}  \left \{
\begin{array}{lll}
\bar{\triangle}_{h} f&=0 &\text{ in }  M,   \\
f&=u  &\text{ on } B.
\end{array}
\right.
\end{equation*}
As the  proof of Theorem 1.1, we show that
\begin{eqnarray}\label{3.1}
  \frac{1}{2}\bar{\triangle}_{h}|\bar{\nabla}f|^{2}&=&|\bar{\nabla}^{2}f|^{2}+\bar{\nabla}\bar{\triangle}_{h}f\bar{\nabla}f+Ric_{h}(\bar{\nabla }f,\bar{\nabla }f)\\
  &\geq& |\bar{\nabla}^{2}f|^{2}+\frac{(\bar{\nabla}h\cdot\bar{\nabla}f)^{2}}{m-n-1}\geq0.\nonumber
\end{eqnarray}
Then by the divergence theorem,  we know that
 $$\int_{B}\frac{\partial }{\partial\nu}(\frac{1}{2}|\bar{\nabla}f|^{2})dA_{h}\geq0.$$
 On the other hand,
\begin{eqnarray}\label{3.2}
  \frac{\partial }{\partial\nu}(\frac{1}{2}|\bar{\nabla}f|^{2})&=&\sum\limits_{i=1}^{n} f_{\nu i}f_{i}-\Pi_{ij}f_{i}f_{j}+f_{\nu\nu}f_{\nu} \\
  &=& \sum\limits_{i=1}^{n} f_{\nu i}u_{i}-\Pi_{ij}u_{i}u_{j}-\triangle_{h}uf_{\nu}-nH_{h}f_{\nu}^{2} \nonumber\\
  &\leq&\sum\limits_{i=1}^{n} f_{\nu i}u_{i}-c|\nabla u|^{2}+\lambda_{1}uf_{\nu}-(m-1)cf_{\nu}^{2}.\nonumber
\end{eqnarray}
Thus we deduce that
\begin{eqnarray}\int_{B}(-\sum\limits_{i=1}^{n} f_{\nu i}u_{i}+c|\nabla u|^{2}-\lambda_{1}uf_{\nu}+(m-1)cf_{\nu}^{2})dA_{h}\leq0.\label{3.3}
\end{eqnarray}
From the equation of the eigenfunction, we see that
\begin{eqnarray*}
  0&\geq&\int_{B}\{-\sum\limits_{i=1}^{n} f_{\nu i}u_{i}+c|\nabla u|^{2}-\lambda_{1}uf_{\nu}+(m-1)cf_{\nu}^{2}\}dA_{h}\\
  &=&\int_{B}\{\triangle_{h} u f_{\nu }+c\lambda_{1}u^{2}-\lambda_{1}uf_{\nu}+(m-1)cf_{\nu}^{2}\}dA_{h}\nonumber \\
 &=&\int_{B}\{c\lambda_{1}u^{2}-2\lambda_{1}uf_{\nu}+(m-1)cf_{\nu}^{2}\}dA_{h} \nonumber\\
  &=&\int_{B}\{(m-1)c(f_{\nu}-\frac{\lambda_{1}}{(m-1)c}u)^{2}+\lambda_{1}(c-\frac{\lambda_{1}}{(m-1)c})u^{2}\}dA_{h}\nonumber\\
  &\geq&\lambda_{1}(c-\frac{\lambda_{1}}{(m-1)c})\int_{B}u^{2}dA_{h}.
\end{eqnarray*}

Thus we have
\begin{equation*}
\lambda_{1}\geq(m-1)c^{2}.
\end{equation*}

If M is isometric to an Euclidean ball of radius $\frac{1}{c}$ and $h$ is a constant, then $m=n+1$, $\lambda_{1}(\triangle_{h})=\lambda_{1}(\triangle)=nc^{2}=(m-1)c^{2}.$ Now we assume that $\lambda_{1}(\triangle_{h})=(m-1)c^{2}.$ In this case, the equality signs in \eqref{3.1},   \eqref{3.2} and  \eqref{3.3} hold. In particular, we see
\begin{equation*}
\bar{\nabla}^{2}f=0, H_{h}=\frac{(m-1)}{n}c, f_{\nu}=\frac{\lambda_{1}}{(m-1)c}u=cu.
\end{equation*}
Since $f$ is the first eigenfunction on B, then $f$ is not a constant. Therefore from $\bar{\nabla}^{2}f=0$, we know $|\bar{\nabla}f|^{2}$ is a constant. By scaling, we  assume $|\bar{\nabla}f|^{2}=1$. Thus  we see that
\begin{eqnarray}
A_{h}=\int_{B}|\bar{\nabla}f|^{2}dA_{h}=\int_{B}(|\nabla u|^{2}+f_{\nu}^{2})dA_{h}=\int_{B}(\lambda_{1}u^{2}+f_{\nu}^{2})dA_{h}=mc^{2}\int_{B}u^{2}dA_{h}.\label{3.4}
\end{eqnarray}
On the other hand, since
\begin{eqnarray*}
\frac{1}{2}\bar{\triangle}_{h}(f^{2})=|\bar{\nabla} f|^{2}+f\bar{\triangle}_{h}f=1,
\end{eqnarray*}
then by the divergence theorem we know
\begin{eqnarray}
V_{h}=\int_{M}\frac{1}{2}\bar{\triangle}_{h}(f^{2})dV_{h}=\int_{B}uf_{\nu}dA_{h}=c\int_{B}u^{2}dA_{h}.\label{3.5}
\end{eqnarray}
Combining \eqref{3.4} and \eqref{3.5}, we obtain

\begin{equation*}
H_{h}=\frac{(m-1)}{n}c=\frac{(m-1)A_{h}}{mnV_{h}}.
\end{equation*}
So from Corollary 1.2 we know that M is isometric to an Euclidean ball and $h$ is a constant. Since  $\lambda_{1}(\triangle)=nc^{2}$, then the radius of M is $\frac{1}{c}$.

{\bf Proof of Theorem 1.7. } Let $u$ be an eigenfunction corresponding to the first nonzero eigenvalue $\lambda_{1}$ of the drifting Laplacian of B:
\begin{equation*}\triangle_{h}u=-\lambda_{1}u.
\end{equation*}
Let $f$ be the smooth solution of the Dirichlet problem:

\begin{equation*}  \left \{
\begin{array}{lll}
\bar{\triangle}_{h} f&=0 &\text{ in }  M,   \\
f&=u  &\text{ on } B.
\end{array}
\right.
\end{equation*}
As the proof of Theorem 1.6, we have that
\begin{eqnarray}\label{3.6}
 \frac{1}{2}\bar{\triangle}_{h}|\bar{\nabla}f|^{2}&=&|\bar{\nabla}^{2}f|^{2}+\bar{\nabla}\bar{\triangle}_{h}f\bar{\nabla}f+Ric_{h}(\bar{\nabla }f,\bar{\nabla }f) \\
 &\geq& |\bar{\nabla}^{2}f|^{2}+\frac{(\bar{\nabla}h\cdot\bar{\nabla}f)^{2}}{m-n-1}\geq0. \nonumber
\end{eqnarray}
Then by the divergence theorem,  we know that $$\int_{B}\frac{\partial }{\partial\nu}(\frac{1}{2}|\bar{\nabla}f|^{2})dA_{h}\geq0.$$
 On the other hand,
\begin{eqnarray}\label{3.7}
 \frac{\partial }{\partial\nu}(\frac{1}{2}|\bar{\nabla}f|^{2})= \sum\limits_{i=1}^{n} f_{\nu i}u_{i}-\Pi_{ij}u_{i}u_{j}-\triangle_{h}uf_{\nu}-nH_{h}f_{\nu}^{2} \\
 \leq\sum\limits_{i=1}^{n} f_{\nu i}u_{i}-c|\nabla u|^{2}+\lambda_{1}uf_{\nu}-\frac{\lambda_{1}}{c}f_{\nu}^{2}.\nonumber
\end{eqnarray}

Thus we obtain that
\begin{eqnarray}\label{3.8}
  0&\geq&\int_{B}\{-\sum\limits_{i=1}^{n} f_{\nu i}u_{i}+c|\nabla u|^{2}-\lambda_{1}uf_{\nu}+\frac{\lambda_{1}}{c}f_{\nu}^{2}\}dA_{h} \\
  &=&\int_{B}\{\triangle_{h} u f_{\nu }+c\lambda_{1}u^{2}-\lambda_{1}uf_{\nu}+\frac{\lambda_{1}}{c}f_{\nu}^{2}\}dA_{h}\nonumber\\
  &=&\int_{B}\{c\lambda_{1}u^{2}-2\lambda_{1}uf_{\nu}+\frac{\lambda_{1}}{c}f_{\nu}^{2}\}dA_{h} \nonumber\\
  &=&\int_{B}\frac{\lambda_{1}}{c}(cu-f_{\nu})^{2}\}dA_{h}\geq0. \nonumber
\end{eqnarray}

Thus  the equalities sign in \eqref{3.6}, \eqref{3.7}  and \eqref{3.8}  hold. In particular, we see
\begin{equation*}
\bar{\nabla}^{2}f=0, H_{h}=\frac{\lambda_{1}}{nc}, f_{\nu}=cu.
\end{equation*}
As the proof of Theorem 1.6, we   assume $|\bar{\nabla}f|^{2}=1$. Thus  we deduce that
\begin{eqnarray}
A_{h}=\int_{B}|\bar{\nabla}f|^{2}dA_{h}=\int_{B}(|\nabla u|^{2}+f_{\nu}^{2})dA_{h}=\int_{B}(\lambda_{1}u^{2}+f_{\nu}^{2})dA_{h}=(\lambda_{1}+c^{2})\int_{B}u^{2}dA_{h}.\label{3.10}
\end{eqnarray}
On the other hand, since
\begin{equation*}
\frac{1}{2}\bar{\triangle}_{h}(f^{2})=|\bar{\nabla} f|^{2}+f\bar{\triangle}_{h}f=1,
\end{equation*}
then by the divergence theorem we know
\begin{eqnarray}
V_{h}=\int_{M}\frac{1}{2}\bar{\triangle}_{h}(f^{2})dV_{h}=\int_{B}uf_{\nu}dA_{h}=c\int_{B}u^{2}dA_{h}.\label{3.11}
\end{eqnarray}
Combining  \eqref{3.10} and  \eqref{3.11}, we obtain

\begin{eqnarray}
\frac{A_{h}}{V_{h}}=\frac{\lambda_{1}+c^{2}}{c}.\label{3.12}
\end{eqnarray}
Notice that $H_{h}\geq\frac{\lambda_{1}}{nc}$, then from Theorem 1.1, we see that
\begin{eqnarray}
\lambda_{1}\leq\frac{(m-1)c}{m}\cdot\frac{A_{h}}{V_{h}}.\label{3.13}
\end{eqnarray}
By \eqref{3.12} and \eqref{3.13}, we obtain that
\begin{equation*}
\lambda_{1}\leq(m-1)c^{2}.
\end{equation*}
Thus we have completed the proof of Theorem 1.7.


\begin{thebibliography}{99}


\bibitem{CMCR}A. Canete, M. Miranda  and  D. Vittone, Some isoperimetric problems in planes with density, J.
Geom. Anal. 20 (2010), No. 2, 243-290.


\bibitem{E} J. F. Escobar, Uniqueness theorems on conformal deformation of metrics, Sobolev inequalities, and an eigenvalue estimate, Comm. Pure Appl. Math., 43 (1990),  857-883.

\bibitem{Espinar}J.M. Espinar, Manifolds with density, applications and gradient Schr\"{o}dinger operators, arXiv preprint arXiv:1209.6162, 2012.


\bibitem{Ho} P. T. Ho, On the structure of $\phi$ stable minimal hypersurfaces in manifolds of nonnegative P-scalar curvature, Math. Ann., 348 (2010), 319-332.


\bibitem{Li}X.D.Li, Liouville theorems for symmetric diffusion
operators on complete Riemannian manifolds, J.Math.Pures Appl.,
84 (2005), 1295-1361.

\bibitem{LW}H.Z. Li and Y. Wei, $f$-minimal surface and manifold with positive m-Bakery-Emery Ricci curvature, arXiv preprint arXiv:1209.0895, 2012.

\bibitem{Lo}J.Lott, some geometric properties of the Bakry-Emery
Ricci tensor, Comment. Math.Helv., 78 (2003), 865-883.

\bibitem{MD}Li Ma and S. H. Du, Extension of Reilly formula with applications to eigenvalue
estimates for drifting Laplacians, C.R.Acad.Sci.Paris, Ser.I, 348 (2010), 1203-1206.

\bibitem{MM}Q. Maurmann and  F. Morgan, Isoperimetric comparison theorems for manifolds with density. Calc.
Var. Partial Differential Equations, 36 (2009), No. 1, 1-5.

\bibitem{Morgan}  F. Morgan, Manifolds with density. Notices Amer. Math. Soc., 52 (2005),
853-858.

\bibitem{Q} Z., Qian, Estimates for weighted volumes and applications.
Quart. J. Math.,   2 (1997)48,  235--242.

\bibitem{Ros}A. Ros, Compact hypersurfaces with constant higher order mean curvatures, Rev. Mat.
Iberoamericana, 3 (1987), 447-453.


\bibitem{Reilly1}R. Reilly, Geometric applications of the solvability of Newmann problems on a Riemannian manifold, Arch. Rat.
Mech. Anal., 75 (1980), 23-30.

\bibitem{R2}R. Reilly, Applications of the Hessian operator in a Riemannian manifold, Indiana Univ.
Math. J., 26 (1977), 459-472.


\bibitem{RCBM} C. Rosales, A. Canete, V. Bayle, and F. Morgan, On the isoperimetric problem in Euclidean space
with density. Calc. Var. Partial Differential Equations, 31 (2008), No.1, 27-46.


\bibitem{Ruan1}Q.H. Ruan, Two rigidity theorems on manifolds with Bakry-Emery Ricci curvature,
Proc. Japan Acad., Ser. A, 85 (2009), 71-74.

\bibitem{Ruan2}Q.H. Ruan, Sharp gradient estimate for positive solutions of the Laplacian with drift, Chin. Ann. Math., Ser.A,  29 (2008), No.1, 107-122.

\bibitem{xia}C. Xia, Rigidity of compact manifolds with boundary and nonnegative Ricci
curvature, Proc. Amer. Math. Soc., 125 (1997), 1801-806.




\end{thebibliography}
\end{document}